\documentclass[12pt]{article}%
\usepackage{amsmath}
\usepackage{amsfonts}
\usepackage{amssymb}
\usepackage{graphicx}%
\providecommand{\U}[1]{\protect\rule{.1in}{.1in}}
\newtheorem{theorem}{Theorem}[section]

\newtheorem{corollary}[theorem]{Corollary}

\newtheorem{lemma}[theorem]{Lemma}

\newtheorem{observation}[theorem]{Observation}
\newtheorem{problem}{Problem}
\newtheorem{proposition}[theorem]{Proposition}
\newtheorem{question}{Question}

\begin{document}

\title{\textbf{A Sharp Upper Bound for the Boundary Independence Broadcast Number of
a Tree}}
\author{C.M. Mynhardt\thanks{Supported by the Natural Sciences and Engineering
Research Council of Canada, PIN 253271.}\\Department of Mathematics and Statistics\\University of Victoria, Victoria, BC, \textsc{Canada}\\{\small kieka@uvic.ca}
\and L. Neilson\\Department of Adult Basic Education\\Vancouver Island University, Nanaimo, BC, \textsc{Canada}\\{\small linda.neilson@viu.ca}}
\maketitle

\begin{abstract}
A broadcast on a nontrivial connected graph $G=(V,E)$ is a function
$f:V\rightarrow\{0,1,\dots,\operatorname{diam}(G)\}$ such that $f(v)\leq e(v)$
(the eccentricity of $v$) for all $v\in V$. The weight of $f$ is $\sigma(f)=%
{\textstyle\sum_{v\in V}}
f(v)$. A vertex $u$ hears $f$ from $v$ if $f(v)>0$ and $d(u,v)\leq f(v)$.

A broadcast $f$ is boundary independent\emph{ }if, for any vertex $w$ that
hears $f$ from vertices $v_{1},...,v_{k},\ k\geq2$, $d(w,v_{i})=f(v_{i})$ for
each $i$. The maximum weight of a boundary independent broadcast on $G$ is
denoted by\emph{ }$\alpha_{\operatorname{bn}}(G)$. We prove a sharp upper
bound on $\alpha_{\operatorname{bn}}(T)$ for a tree $T$ in terms of its order
and number of branch vertices of a certain type.

\end{abstract}

\noindent\textbf{Keywords:\hspace{0.1in}}broadcast domination; broadcast
independence, hearing independence; boundary independence

\noindent\textbf{AMS Subject Classification Number 2010:\hspace{0.1in}}05C69

\section{Introduction}

Boundary independent broadcasts on graphs were defined by Mynhardt and Neilson
\cite{MN} and Neilson \cite{LindaD} as an alternative to the concept of
independent broadcasts, here called hearing independent broadcasts, as defined
by Erwin \cite{Ethesis}. Both concepts are generalizations of independent sets
in graphs. The former focusses on the fact that no edge is covered by more
than one vertex of an independent set $X$, and the latter on the property that
no vertex in $X$ belongs to the neighbourhood of (or hears) another vertex in
$X$. As a result, in a boundary independent broadcast, any edge is covered by
at most one broadcasting vertex, and a vertex belongs to the broadcasting
neighbourhoods of two or more broadcasting vertices only if it belongs to the
boundaries of all such vertices, whereas in a hearing independent broadcast it
is only required that no broadcasting vertex belongs to the broadcasting
neighbourhood of another such vertex.

Our purpose is to prove a sharp upper bound for the boundary independence
number of a tree in terms of its order and the number of branch vertices with
certain properties. After giving the necessary definitions in Sections
\ref{Sec_broadcast_Defs} and \ref{SecTreeDefs}, we state the bound in Section
\ref{SecMainThm} (Theorem~\ref{ThmMain}). The proof is given in Section
\ref{SecMainThmPf}. Sections \ref{Sec_known} and \ref{SecLemmas} contain
previous results and lemmas required for the proof of the bound. For a class
of trees that meet the bound we briefly consider boundary independent
broadcasts on caterpillars in Section \ref{Sec_Cat}, and conclude by listing
open problems in Section \ref{Sec_Open}. For undefined concepts we refer the
reader to \cite{CLZ}.

\subsection{Broadcast definitions}

\label{Sec_broadcast_Defs}The study of broadcast domination was initiated by
Erwin \cite{Ethesis, Epaper}. A \emph{broadcast} on a nontrivial connected
graph $G=(V,E)$ is a function $f:V\rightarrow\{0,1,\dots,\operatorname{diam}%
(G)\}$ such that $f(v)\leq e(v)$ (the eccentricity of $v$) for all $v\in V$.
If $G$ is disconnected, we define a broadcast on $G$ as the union of
broadcasts on its components. The \emph{weight} of $f$ is $\sigma
(f)=\sum_{v\in V}f(v)$. Define $V_{f}^{+}=\{v\in V:f(v)>0\}$ and partition
$V_{f}^{+}$ into the two sets $V_{f}^{1}=\{v\in V:f(v)=1\}$ and $V_{f}%
^{++}=V_{f}^{+}-V_{f}^{1}$. A vertex in $V_{f}^{+}$ is called a
\emph{broadcasting vertex}. A vertex $u$ \emph{hears} $f$ from $v\in V_{f}%
^{+}$, and $v$ $f$-\emph{dominates} $u$, if the distance $d(u,v)\leq f(v)$. If
$d(u,v)<f(v)$, we also say that say that $v$ \emph{overdominates }$u$. Denote
the set of all vertices that do not hear $f$ by $U_{f}$. A broadcast $f$ is
\emph{dominating} if $U_{f}=\varnothing$. If $f$ is a broadcast such that
every vertex $x$ that hears more than one broadcasting vertex also satisfies
$d(x,u)\geq f(u)$ for all $u\in V_{f}^{+}$, we say that the \emph{broadcast
only overlaps in boundaries}. If $uv\in E(G)$ and $u,v\in N_{f}(x)$ for some
$x\in V_{f}^{+}$ such that at least one of $u$ and $v$ does not belong to
$B_{f}(x)$, we say that the edge $uv$ is \emph{covered} in $f$, or
$f$-\emph{covered}, by $x$. If $uv$ is not covered by any $x\in V_{f}^{+}$, we
say that $uv$ is \emph{uncovered by~}$f$ or $f$-\emph{uncovered}. We denote
the set of $f$-uncovered edges by $U_{f}^{E}$.

If $f$ and $g$ are broadcasts on $G$ such that $g(v)\leq f(v)$ for each $v\in
V$, we write $g\leq f$. If in addition $g(v)<f(v)$ for at least one $v\in V$,
we write $g<f$. A dominating broadcast $f$ on $G$ is a \emph{minimal
dominating broadcast} if no broadcast $g<f$ is dominating. The \emph{upper
broadcast number }of $G$ is
\[
\Gamma_{b}(G)=\max\left\{  \sigma(f):f\text{ is a minimal dominating broadcast
of }G\right\}  ,
\]
and a dominating broadcast $f$ of $G$ such that $\sigma(f)=\Gamma_{b}(G)$ is
called a $\Gamma_{b}$-\emph{broadcast}. First defined by Erwin \cite{Ethesis},
the upper broadcast number was also studied by, for example, Ahmadi, Fricke,
Schroeder, Hedetniemi and Laskar \cite{Ahmadi}, Bouchemakh and Fergani
\cite{BF}, Bouchouika, Bouchemakh and Sopena \cite{BBS}, Dunbar, Erwin,
Haynes, Hedetniemi and Hedetniemi \cite{DEHHH}, and Mynhardt and Roux
\cite{MR}.

We denote the independence number of $G$ by $\alpha(G)$; an independent set of
$G$ of cardinality $\alpha(G)$ is called an $\alpha$-\emph{set}. If $f$ is
characteristic function of an independent set of $G$, then no vertex in
$V_{f}^{+}$ hears $f$ from any other vertex. To generalize the concept of
independent sets, Erwin \cite{Ethesis} defined a broadcast $f$ to be
\emph{independent}, or, for our purposes, \emph{hearing independent},
abbreviated to \emph{h-independent}, if no vertex $u\in V_{f}^{+}$ hears $f$
from any other vertex $v\in V_{f}^{+}$; that is, broadcasting vertices only
hear themselves. The maximum weight of an h-independent broadcast is the
\emph{h-independent broadcast number}, which we denote by $\alpha_{h}(G)$;
such a broadcast is called an $\alpha_{h}$-\emph{broadcast}. This version of
broadcast independence was also considered by, among others, Ahmane,
Bouchemakh and Sopena \cite{ABS, ABS2}, Bessy and Rautenbach \cite{BR, BR2},
Bouchemakh and Zemir \cite{Bouch}, Bouchouika et al.~\cite{BBS} and Dunbar et
al.~\cite{DEHHH}. For a survey of broadcasts in graphs, see the chapter by
Henning, MacGillivray and Yang \cite{HMY}.

For a broadcast $f$ on a graph $G$ and $v\in V_{f}^{+}$, we define the%
\[
\left.
\begin{tabular}
[c]{l}%
$f$-\emph{neighbourhood}\\
$f$-\emph{boundary}\\
$f$-\emph{private neighbourhood}\\%
\begin{tabular}
[c]{l}%
\
\end{tabular}
\\
$f$-\emph{private boundary}\\%
\begin{tabular}
[c]{l}%
\
\end{tabular}
\end{tabular}
\ \ \ \right\}  \ \text{of }v\text{ by }\left\{  \text{%
\begin{tabular}
[c]{rll}%
$N_{f}(v)$ & $=$ & $\{u\in V:d(u,v)\leq f(v)\}$\\
$B_{f}(v)$ & $=$ & $\{u\in V:d(u,v)=f(v)\}$\\
$\operatorname{PN}_{f}(v)$ & $=$ & $\{u\in N_{f}(v):u\notin N_{f}(w)$ for
all\\
&  & $\ w\in V_{f}^{+}-\{v\}\}$\\
$\operatorname{PB}_{f}(v)$ & $=$ & $\{u\in N_{f}(v):u$ is not dominated by\\
&  & $\ (f-\{(v,f(v)\})\cup\{(v,f(v)-1)\}$.
\end{tabular}
}\right.
\]

If $v\in V_{f}^{1}$ and $v$ does not hear $f$ from any vertex $u\in V_{f}%
^{+}-\{v\}$, then $v\in\operatorname{PB}_{f}(v)$, and if $v\in V_{f}^{++}$,
then $\operatorname{PB}_{f}(v)=B_{f}(v)\cap\operatorname{PN}_{f}(v)$. Also
note that $f$ is a broadcast that overlaps only in boundaries if and only if
$N_{f}(u)\cap N_{f}(v)\subseteq B_{f}(u)\cap B_{f}(v)$ for all distinct
$u,v\in V_{f}^{+}$.

The characteristic function of an independent set also has the feature that it
only overlaps in boundaries. To generalize this property, we define a
broadcast to be \emph{boundary independent}, abbreviated to
\emph{bn-independent}, if it overlaps only in boundaries. The maximum weight
of a bn-independent broadcast on $G$ is the \emph{boundary independence
number} $\alpha_{\operatorname{bn}}(G)$; such a broadcast is called an
$\alpha_{\operatorname{bn}}(G)$\emph{-broadcast}, often abbreviated to
$\alpha_{\operatorname{bn}}$\emph{-broadcast} if the graph $G$ is clear. The
respective definitions imply that $\alpha_{\operatorname{bn}}(G)\leq\alpha
_{h}(G)$ for all graphs $G$. Boundary independent broadcasts were introduced
by Neilson \cite{LindaD} and Mynhardt and Neilson \cite{MN}, and also studied
in \cite{MM, MN2}. For example, it was shown in \cite{MN} that $\alpha
_{h}(G)/\alpha_{\operatorname{bn}}(G)<2$ for all graphs $G$, and the bound is
asymptotically best possible. In \cite{MN2} it was shown that $\alpha
_{\operatorname{bn}}$ and $\Gamma_{b}$ are not comparable, and that
$\alpha_{\operatorname{bn}}(G)/\Gamma_{b}(G)<2$ for all graphs $G$, while
$\Gamma_{b}(G)/\alpha_{\operatorname{bn}}(G)$ is unbounded.

\subsection{Definitions for trees}

\label{SecTreeDefs}
\begin{figure}[ptb]%
\centering
\includegraphics[
height=3.0217in,
width=5.6265in
]%
{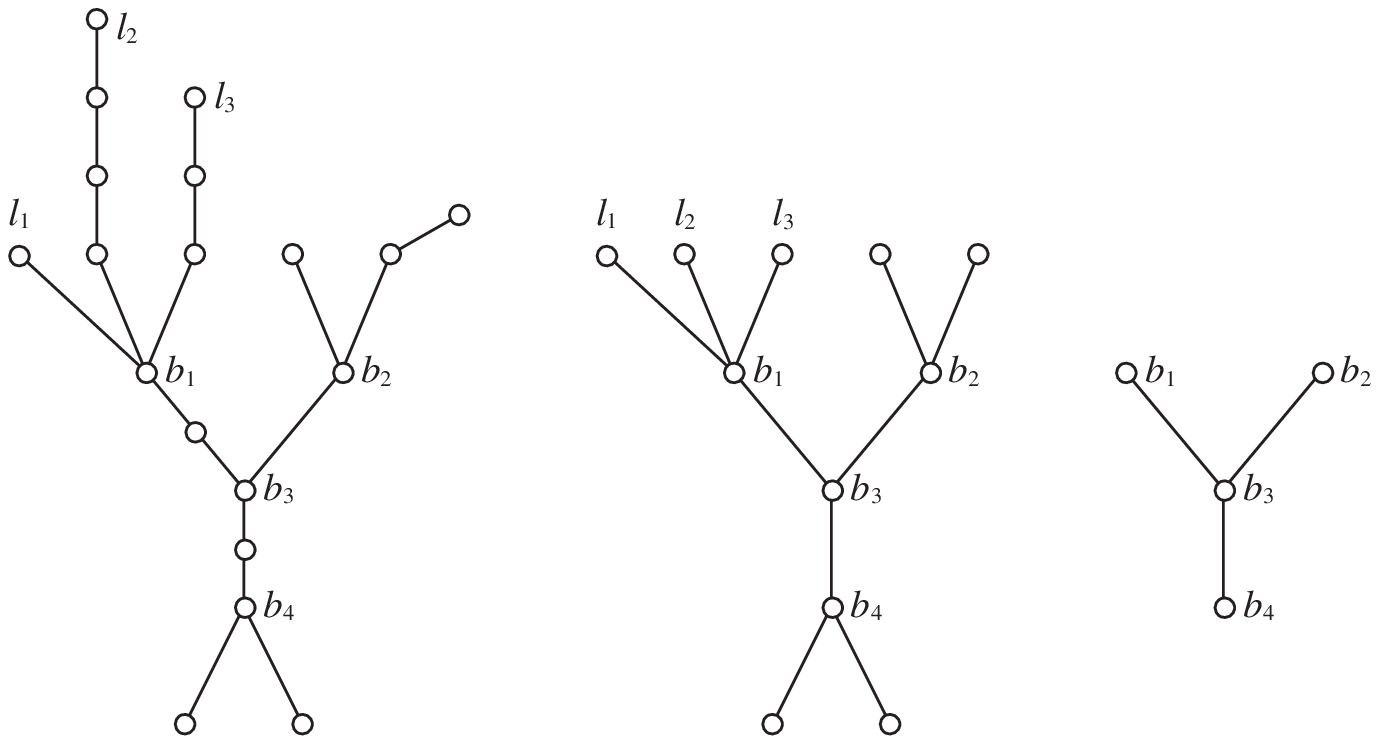}%
\caption{A tree $T$ (left), its branch-leaf representation $\mathcal{BL}(T)$
(middle) and its branch representation $\mathcal{B}(T)$ (right). The branch
set of $T$ is $B(T)=\{b_{1},b_{2},b_{3},b_{4}\}$, $R(T)=\{b_{3}\}$,
$B_{\operatorname{end}}(T)=\{b_{1},b_{2},b_{4}\}$, $b(T)=4$ and $\rho(T)=1$.
The branch vertex $b_{1}$ has leaf set $L(b_{1})=\{l_{1},l_{2},l_{3}\}$. }%
\label{Fig6.1}%
\end{figure}

The statement of our main result requires some definitions
of concepts pertaining to trees. A vertex of a tree $T$ of degree $3$ or more
is called a \emph{branch vertex}. We denote the set of leaves of $T$ by
$L(T)$, the set of branch vertices by $B(T)$ and the set of vertices of degree
$2$ by $W(T)$. The unique neighbour of a leaf is called a \emph{stem}. The
next few concepts are illustrated in Figure \ref{Fig6.1}. The
\emph{branch-leaf representation }$\mathcal{BL}(T)$ of $T$ is the tree
obtained by suppressing all vertices $v$ with $\deg(v)=2$, and the
\emph{branch representation }$\mathcal{B}(T)$ of a tree $T$ with at least one
branch vertex is obtained by deleting all leaves of $\mathcal{BL}(T)$. Thus,
$V(\mathcal{B}(T))=B(T)$, and two vertices $b_{1},b_{2}\in B(T)$ are adjacent
in $\mathcal{B}(T)$ if and only if the $b_{1}-b_{2}$ path in $T$ contains no
other branch vertices. We denote $|B(T)|$ by $b(T)$.%

An \emph{endpath} in a tree is a path ending in a leaf and having all internal
vertices (if any) of degree $2$. If there exists a $v-l$ endpath, where $v\in
B(T)$ and $l\in L(T)$, then $v$ and $l$ are adjacent \underline{in
$\mathcal{BL}(T)$}; we also say that $l$ belongs to $L(v)$, the \emph{leaf
set} of $v$, and we refer to $l$ as \emph{a leaf of} $v$ (even though $l$ is
not necessarily adjacent to $v$ in $T$). Since $\mathcal{BL}(T)$ is unique, we
can talk about $L(v)$ for any branch vertex $v$, where the reference to
$\mathcal{BL}(T)$ is implied but not specifically mentioned. Let $R(T)$ be the
set of all branch vertices $w$ of $T$ such that $|L(w)|\leq1$ and define
$\rho(T)=|R(T)|$. Equivalently, $\rho(T)$ is the number of branch vertices of
$T$ with at most one leaf, that is, the branch vertices which belong to at
most one endpath. A branch vertex of $T$ of degree $1$ in $\mathcal{B}(T)$ is
called an \emph{endbranch vertex}. Denote the set of endbranch vertices by
$B_{\operatorname{end}}(T)$. Since $\mathcal{B}(T)$ is a tree, and any tree of
order at least $2$ has at least two leaves, $|B_{\operatorname{end}}(T)|\geq2$
for every tree $T$ with $b(T)\geq2$. If $w\in B_{\operatorname{end}}(T)$, then
exactly one edge incident with $w$ does not lie on an endpath of $T$. Hence
$|L(w)|\geq2$ and $B_{\operatorname{end}}(T)\cap R(T)=\varnothing$.

We define subsets $B_{i}(T)$ and $B_{\geq i}(T)$ of $B(T)$ by
\[
B_{i}(T)=\{v\in B(T):|L(v)|=i\}\text{ and }B_{\geq i}(T)=\{v\in
B(T):|L(v)|\geq i\}.
\]
Clearly, $B_{0}(T)\cup B_{1}(T)\cup B_{\geq2}(T)$ is a partition of $B(T)$
while $B_{0}(T)\cup B_{1}(T)$ is a partition of $R(T)$. We also partition the
set $W(T)$ of vertices of degree $2$ into two subsets, $W_{\operatorname{ext}%
}(T)$ for the \emph{external vertices of degree }$2$, and
$W_{\operatorname{int}}(T)$ for the \emph{internal vertices of degree }$2$, as
follows:%
\[
W_{\operatorname{ext}}(T)=\{u\in W(T):u\ \text{lies\ on\ an\ endpath}%
\}\text{\ and\ }W_{\operatorname{int}}(T)=W(T)-W_{\operatorname{ext}}(T).
\]
The subgraph of $T$ induced by $B_{0}(T)\cup B_{1}(T)\cup
W_{\operatorname{int}}(T)$ is called the \emph{interior subgraph} of $T$,
denoted by $\operatorname{Int}(T)$. Note that $\operatorname{Int}(T)$ is
acyclic but not necessarily connected.

\subsection{Statement of main theorem}

\label{SecMainThm}
\begin{figure}[ph]%
\centering
\includegraphics[
height=0.7169in,
width=1.9268in
]%
{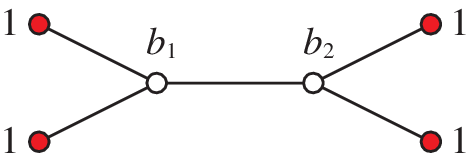}%
\caption{An $\alpha_{\operatorname{bn}}(T)$-broadcast, where the branch
vertices are covered by leaves and the edge between them is uncovered;
$\alpha_{\operatorname{bn}}(T)=n-b(T)=4$.}%
\label{Fig_6.6}%
\end{figure}

\begin{figure}[ph]%
\centering
\includegraphics[
height=0.7169in,
width=1.8161in
]%
{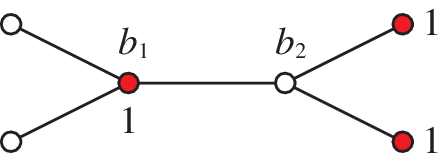}%
\caption{A bn-independent broadcast $f$ in which the edge between two branch
vertices is covered by the branch vertex $b_{1}$, which covers $f(b_{1})+2$
edges. The broadcast is not an $\alpha_{\operatorname{bn}}(T)$-broadcast and
$\sigma(f)=3$.}%
\label{Fig_6.7}%
\end{figure}

\begin{figure}[ph]%
\centering
\includegraphics[
height=0.7152in,
width=2.0911in
]%
{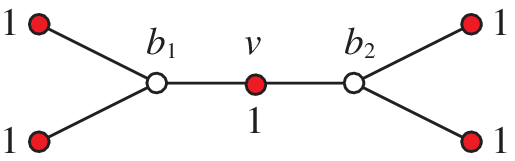}%
\caption{An $\alpha_{\operatorname{bn}}(T)$-broadcast, where the edges between
the branch vertices are covered by a vertex $v\in W(T)$, but $f(v)=1$ and $v$
covers $f(v)+1=2$ edges; $\alpha_{\operatorname{bn}}(T)=n-b(T)=5$.}%
\label{Fig_6.8}%
\end{figure}

\begin{figure}[ph]%
\centering
\includegraphics[
height=0.7031in,
width=2.2459in
]%
{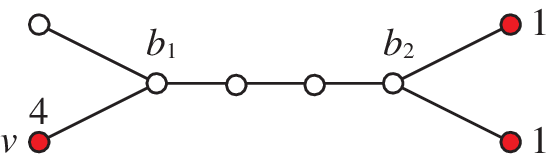}%
\caption{An $\alpha_{\operatorname{bn}}(T)$-broadcast, where the edges between
branch vertices are covered by a leaf $v$; $\alpha_{\operatorname{bn}%
}(T)=n-b(T)=6$.}%
\label{Fig_6.9}%
\end{figure}

It is clear from the definition of bn-independence that if
$f$ is a bn-independent broadcast on an $n$-vertex graph $G$, then the set of
uncovered edges together with the sets of edges covered by $N_{f}(v)$ for all
$v\in V_{f}^{+}$ is a partition of $E(G)$. Thus the upper bound $\alpha
_{\operatorname{bn}}(G)\leq n-1$ is achieved when all edges are covered and
each $v\in V_{f}^{+}$ covers exactly $f(v)$ edges. However, if a tree has two
or more branch vertices then this upper bound is not achievable. We can see
this by considering how the broadcast covers the edges between two branch
vertices. There are four possibilities: the branch vertices are covered by
leaves and the edge between them is uncovered (Figure \ref{Fig_6.6}), the edge
between two branch vertices is covered by a branch vertex (Figure
\ref{Fig_6.7}), the edges between branch vertices $u$ and $w$ are covered by
internal vertices $v_{1},...,v_{k}$ on the $u-w$ path, where $f(v_{i})=1$ and
each $v_{i}$ covers two edges (Figure \ref{Fig_6.8}), and the edges between
branch vertices are covered by a leaf (Figure \ref{Fig_6.9}).
%


\begin{figure}[ptb]%
\centering
\includegraphics[
height=1.5056in,
width=1.4477in
]%
{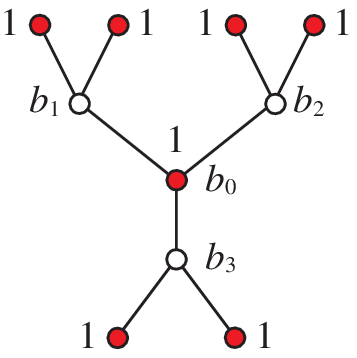}%
\caption{An $\alpha_{\operatorname{bn}}(T)$-broadcast on a tree $T$ in which
$L(b_{0})=\varnothing$, $R(T)=\{b_{0}\}$ and $\alpha_{\operatorname{bn}%
}(T)=7=n-b(T)+\rho(T)>n-b(T)$.}%
\label{Fig_6.10}%
\end{figure}

\begin{figure}[ptb]%
\centering
\includegraphics[
height=1.1917in,
width=2.8746in
]%
{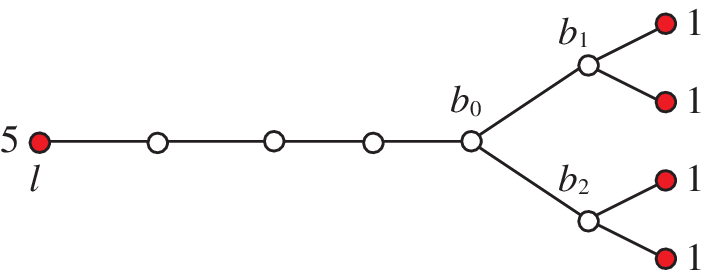}%
\caption{An $\alpha_{\operatorname{bn}}(T)$-broadcast on a tree $T$ in which
$L(b_{0})=\{l\}$, $R(T)=\{b_{0}\}$ and $\alpha_{\operatorname{bn}%
}(T)=9=n-b(T)+\rho(T)>n-b(T)$.}%
\label{Fig_6.11}%
\end{figure}

It may seem that, after the first branch vertex, every additional branch
vertex reduces the value of $\alpha_{\operatorname{bn}}(T)$ relative to $n$.
However, this is not always the case, as is shown by the trees in Figures
\ref{Fig_6.10} and \ref{Fig_6.11}. In Figure \ref{Fig_6.10}, $\mathcal{B}%
(T)\cong K_{1,3}$ and $L(b_{0})=\varnothing$. If $f$ is a bn-independent
broadcast on $T$ such that the central branch vertex $b_{0}$ of the star is a
broadcasting vertex with $f(b_{0})=1$, then $\deg(b_{0})$ edges are covered by
$b_{0}$ while $\deg(b_{0})+1$ branch vertices are dominated. Note that
$R(T)=\{b_{0}\}$ and $\rho(T)=1$. In Figure \ref{Fig_6.11} we see an
$\alpha_{\operatorname{bn}}$-broadcast on a tree $T$ for which $\mathcal{B}%
(T)\cong K_{1,2}$. Note that $L(b_{0})=\{l\}$ and $f(l)=d(l,b_{0})+1$. Also,
$l$ overdominates $b_{0}$ by $1$ and $N_{f}(l)$ covers $f(l)+1$ edges. In each
case, $\alpha_{\operatorname{bn}}(T)>n-b(T)$. This observation was the
motivation for the definition of $R(T)$ and $\rho(T)$.%

We now state our main result. The bound is sharp for generalized spiders (see
Proposition \ref{Prop=}) and for some caterpillars (Corollary \ref{Cor_Cat}).

\begin{theorem}
\label{ThmMain}For any tree $T$ of order $n$, $\alpha_{\operatorname{bn}%
}(T)\leq n-b(T)+\rho(T).$
\end{theorem}

\section{Known results}

\label{Sec_known}In this section we present known results that will be used
later on. It is often useful to know when a bn-independent broadcast $f$ is
\emph{maximal bn-independent}, that is, there does not exist a bn-independent
broadcast $g$ such that $f<g$.

\begin{proposition}
\label{prop-max-bn}

\begin{enumerate}
\item[$(i)$] \emph{\cite{MN}}\hspace{0.1in}A bn-independent broadcast $f$ on a
graph $G$ is maximal bn-indepen-dent if and only if it is dominating and
either $V_{f}^{+}=\{v\}$ or $B_{f}(v)-\operatorname{PB}_{f}(v)\neq\varnothing$
for each $v\in V_{f}^{+}$.

\item[$(ii)$] \emph{\cite{MM}\hspace{0.1in}}Let $f$ be a bn-independent
broadcast on a connected graph $G$ such that $|V_{f}^{+}|\geq2$. Then $f$ is
maximal bn-independent if and only if each component of $G-U_{f}^{E}$ contains
at least two broadcasting vertices.
\end{enumerate}
\end{proposition}

Suppose $f$ is a bn-independent broadcast on $G$ and an edge $uv$ of $G$ is
covered by vertices $x,y\in V_{f}^{+}$. By the definition of covered,
$\{u,v\}\nsubseteq B_{f}(x)$ and $\{u,v\}\subseteq N_{f}(x)\cap N_{f}(y)$.
This violates the bn-independence of $f$. Hence we have the following observation.

\begin{observation}
\label{edge-covered}If $f$ is a bn-independent broadcast on a graph $G$, then
each edge of $G$ is covered by at most one vertex in $V_{f}^{+}$.
\end{observation}

For $k\geq3$ and $n_{i}\geq1,\ i\in\{1,...,k\}$, the \emph{(generalized)
spider} $\operatorname{Sp}(n_{1},...,n_{k})$ is the tree which has exactly one
vertex $b$, called the \emph{head}, with $\deg(b)=k$, and for which the $k$
components of $\operatorname{Sp}(n_{1},...,n_{k})-b$ are paths of lengths
$n_{1}-1,...,n_{k}-1$, respectively. The \emph{legs }$L_{1},...,L_{k}$ of the
spider are the paths from $b$ to the leaves. Let $t_{i}$ be the leaf of
$L_{i},\ i=1,...,k$. If $n_{i}=r$ for each $i$, we write $\operatorname{Sp}%
(n_{1},...,n_{k})=\operatorname{Sp}(r^{k})$.

The following bound on $\alpha_{\operatorname{bn}}(G)$ was proved in \cite{MN,
LindaD}.

\begin{proposition}
\label{Prop=}\emph{\cite{MN, LindaD}}\hspace{0.1in}For any connected graph $G$
of order $n$ and any spanning tree $T$ of $G$, $\alpha_{\operatorname{bn}%
}(G)\leq\alpha_{\operatorname{bn}}(T)\leq n-1$. Moreover, $\alpha
_{\operatorname{bn}}(G)=n-1$ if and only if $G$ is a path or a generalized spider.
\end{proposition}

Note that if $T$ is a tree such that $b(T)\leq2$, then $\rho(T)=0$. If
$b(T)=0$, then $T$ is a path, if $\rho(T)=1$, then $T$ is a generalized
spider, and if $\rho(T)=2$, then, by Proposition \ref{Prop=}, $\alpha
_{\operatorname{bn}}(T)\leq n-2$. This shows that Theorem \ref{ThmMain} holds
for $T$:

\begin{observation}
\label{Ob_Base}If $T$ is a tree such that $b(T)\leq2$, then $\alpha
_{\operatorname{bn}}(T)\leq n-b(T)+\rho(T)$.
\end{observation}

Finally, if $P$ is a diametrical path of a tree $T$, then each branch vertex
of $T$ is incident with at least one unique edge not on $P$. By counting
edges, we obtain the following observation.

\begin{observation}
\label{Obs_diam}For any tree $T$ of order $n$, $\operatorname{diam}(T)\leq
n-b(T)-1$.
\end{observation}

\section{Lemmas}

\label{SecLemmas}This section contains a number of lemmas required for the
proof of our main result, Theorem \ref{ThmMain}. Although they were developed
with this specific purpose in mind, they are also useful for finding
$\alpha_{\operatorname{bn}}$-broadcasts on specific trees. We begin by showing
that if $f$ is an $\alpha_{\operatorname{bn}}(T)$-broadcast, then any leaf of
$T$ hears $f$ from a leaf (either itself or another leaf).

\begin{lemma}
\label{Lem_hear1}For any tree $T$ and any $\alpha_{\operatorname{bn}}%
(T)$-broadcast $f$, no leaf of $T$ hears $f$ from any non-leaf vertex.
\end{lemma}

\noindent\textbf{Proof.\hspace{0.1in}}It follows from the proof of Proposition
\ref{Prop=} (and is easy to see directly) that the statement is true for
paths, hence assume $T$ has at least one branch vertex. Suppose there exists
an $\alpha_{\operatorname{bn}}(T)$-broadcast $f$ such that a leaf hears a
vertex $v$ with $\deg(v)\geq2$. Either $v$ overdominates all leaves in
$N_{f}(v)$ or there is at least one leaf $u$ such that $d(u,v)=f(v)$. We
consider the two cases separately.

\noindent\textbf{Case 1:\hspace{0.1in}}There exists a leaf $u$ such that
$d(u,v)=f(v)$. If the $u-v$ path contains a branch vertex, let $w$ be the
branch vertex on this path nearest to $u$; if no such vertex exists, then let
$w=v$. Since $v$ is not a leaf, $\deg(w)\geq2$ in either case. Define the
broadcast $g_{1}$ by%
\[
g_{1}(x)=\left\{
\begin{tabular}
[c]{ll}%
$f(v)-d(u,w)$ & if $x=v$\\
$d(u,w)$ & if $x=u$\\
$f(x)$ & otherwise.
\end{tabular}
\ \right.
\]
Notice that $N_{g_{1}}(u)\cup N_{g_{1}}(v)\subseteq N_{f}(v)$ and either
$g_{1}(v)=0$ or $N_{g_{1}}(u)\cap N_{g_{1}}(v)=\{w\}$. Hence $g_{1}$ is
bn-independent and $\sigma(g_{1})=\sigma(f)$. If $f(v)>0$, then $\deg(w)\geq3$
and $w$ has a $g_{1}$-undominated neighbour. If $f(v)=0$, then $v=w$ and,
since $v$ is not a leaf, it has a $g_{1}$-undominated neighbour. In either
case (the contrapositive of) Proposition \ref{prop-max-bn}$(i)$ implies that
$g_{1}$ can be extended to produce a bn-independent broadcast of larger weight
than $\sigma(f)$, contradicting the maximality of $f$.$\smallskip$

\noindent\textbf{Case 2:\hspace{0.1in}}$v$ overdominates all leaves in
$N_{f}(v)$. There are two subcases.

\noindent\textbf{Case 2a:\hspace{0.1in}}There exist $a,b\in B_{f}(v)$ such
that $P=(a=v_{0},...,v_{2f(v)}=b)$ is a path of length $2f(v)$ in
$T[N_{f}(v)]$ containing $v$. Since $v$ overdominates all leaves in $N_{f}%
(v)$, $a$ and $b$ are not leaves; clearly, the internal vertices of $P$ are
not leaves either. Define the broadcast $g_{2}$ by%
\[
g_{2}(x)=\left\{
\begin{tabular}
[c]{ll}%
$0$ & if $x=v_{i}$ and $i$ is even\\
$1$ & if $x=v_{i}$ and $i$ is odd\\
$f(x)$ & otherwise.
\end{tabular}
\ \right.
\]
Observe that $g_{2}(v_{0})=g_{2}(v_{2f(v)})=0$, $\bigcup_{i=0}^{2f(v)}%
N_{g_{2}}(v_{i})\subseteq N_{f}(v)$ and no adjacent vertices are broadcasting
vertices. Hence $g_{2}$ is bn-independent. Since $\sigma(g_{2})=\sigma
(f)=\alpha_{\operatorname{bn}}(T)$, $g_{2}$ is maximal bn-independent. By
Proposition \ref{prop-max-bn}$(i)$, $g_{2}$ is dominating and every leaf in
$N_{f}(v)$ hears $g_{2}$ from some $v_{i},\ i=0,...,v_{2f(v)}$. Since no
$v_{i}$ is a leaf and $g_{2}(v_{i})\in\{0,1\}$, $g_{2}$ satisfies the
conditions of Case 1 and we obtain a contradiction as before.$\smallskip$

\noindent\textbf{Case 2b:\hspace{0.1in}}$v$ overdominates all leaves in
$N_{f}(v)$ and there is no path of length $2f(v)$ in $T[N_{f}(v)]$. Then there
is exactly one edge $e$ incident with $v$ which is on all $v-w$ paths for
$w\in B_{f}(v)$. Let $u$ be a neighbour of $v$ that is not incident with $e$;
$u$ exists because $\deg(v)\geq2$. Define the broadcast $g_{3}$ by%
\[
g_{3}(x)=\left\{
\begin{tabular}
[c]{ll}%
$f(v)+1$ & if $x=u$\\
$0$ & if $x=v$\\
$f(x)$ & otherwise.
\end{tabular}
\right.
\]
Since all vertices in the subtree of $T-uv$ containing $u$ are at distance
less than $f(v)$ from $v$, $N_{g_{3}}(u)\subseteq N_{f}(v)$, hence $g_{3}$ is
bn-independent and $\sigma(g_{3})>\sigma(f)$, a contradiction.~$\blacksquare$

\bigskip

In our next lemma we show that any tree has an $\alpha_{\operatorname{bn}}%
$-broadcast $f$ such that any non-leaf vertex in $V_{f}^{+}$ belongs to
$V_{f}^{1}$.

\begin{lemma}
\label{Lem_nonleaf1}For any tree $T$ there exists an $\alpha
_{\operatorname{bn}}(T)$-broadcast $f$ such that $f(v)=1$ whenever $v\in
V_{f}^{+}-L(T)$.
\end{lemma}

\noindent\textbf{Proof.\hspace{0.1in}}Suppose there exists a tree\textbf{ }$T$
for which the statement of the lemma is false. Among all $\alpha
_{\operatorname{bn}}(T)$-broadcasts, let $f$ be one such that the non-leaf
vertices in $V_{f}^{++}$ is minimum. Let $v\in V_{f}^{++}$ be a non-leaf
vertex. By Lemma \ref{Lem_hear1}, $v$ does not broadcast to a leaf. Therefore
there exist vertices $a,b\in B_{f}(v)$ such that the $a-b$ path $P=(a=v_{0}%
,...,v_{2f(v)}=b)$ contains $v$. Now the $\alpha_{\operatorname{bn}}%
(T)$-broadcast $g_{2}$ defined as in Case 2a of the proof of Lemma
\ref{Lem_hear1} has fewer non-leaf vertices in $V_{g_{2}}^{++}$ than there are
in $V_{f}^{++}$, contradicting the choice of $f$.~$\blacksquare$

\bigskip

We make one more observation about the structure of $\alpha_{\operatorname{bn}%
}$-broadcasts on trees.

\begin{lemma}
\label{Lem_PB_empty}Let $f$ be an $\alpha_{\operatorname{bn}}$-broadcast on a
tree $T$ such that $|V_{f}^{1}|$ is maximum. Then $\operatorname{PB}%
_{f}(v)=\varnothing$ for all $v\in V_{f}^{++}$.
\end{lemma}

\noindent\textbf{Proof.\hspace{0.1in}}Suppose, contrary to the statement, that
there exist vertices $v\in V_{f}^{++}$ and $u\in\operatorname{PB}_{f}(v)$.
Define the broadcast $g$ by $g(x)=f(x)-1$ if $x=v$, $g(x)=1$ if $x=u$, and
$g(x)=f(x)$ otherwise. Since $u\in\operatorname{PB}_{f}(v)$, $u$ does not hear
$g$ from $v$ or any other vertex in $V_{f}^{+}$, hence $g$ is bn-independent.
Since $\sigma(g)=\sigma(f)$, $g$ is an $\alpha_{\operatorname{bn}}%
(T)$-broadcast. But $|V_{g}^{1}|>|V_{f}^{1}|$, contradicting the choice of
$f$.~$\blacksquare$

\bigskip

In the next two lemmas we develop results regarding the role of broadcasting
leaves in $\alpha_{\operatorname{bn}}(T)$-broadcasts. The first one states
essentially that each endpath in $T$ contains at most one broadcasting vertex
of an $\alpha_{\operatorname{bn}}(T)$-broadcast, namely its leaf.

\begin{lemma}
\label{Lem_leaf_branch}Let $f$ be an $\alpha_{\operatorname{bn}}$-broadcast on
a tree $T$ such that a leaf $l$ $f$-dominates a branch vertex $w$. If
$l^{\prime}$ is a leaf in $L(w)$ that does not hear $f$ from $l$, then
$l^{\prime}\in V_{f}^{+}$, the $l^{\prime}-w$ path $Q$ contains a vertex $b\in
B_{f}(l)$, and $f(l^{\prime})=d(b,l^{\prime})$.
\end{lemma}

\noindent\textbf{Proof.\hspace{0.1in}}Let $e$ be the edge incident with $w$ on
the $l-w$ path that does not belong to $Q$; note that $e$ is $f$-covered by
$l$. Since $l^{\prime}\in L(w)$, each internal vertex of $Q$ has degree $2$ in
$T$. Since $l$ dominates $w$ but not $l^{\prime}$, the path $Q$ contains a
vertex $b\in B_{f}(l)$. (Possibly, $b=w$.) Since $f$ is dominating and leaves
only hear leaves (Lemma \ref{Lem_hear1}), some leaf $t\neq l$ broadcasts to
$l^{\prime}$. If $t$ does not belong to $Q$, then $t$ overdominates $w$. But
then $t$ also $f$-covers $e$, which is impossible. Hence $t\in V(Q)$ and so
$t=l^{\prime}$. Let $a$ be the broadcast vertex on the $l^{\prime}-b$ subpath
$R$ of $Q$ nearest to $b$. Since $b\in B_{f}(l)$, $a\neq b$.

Suppose $a\neq l^{\prime}$. Let $a_{1},...,a_{k}$ be all the broadcast
vertices on $R$ strictly between $l^{\prime}$ and $b$. Since $f$ is
bn-independent, $a_{1},...,a_{k}$ cover exactly $2\sum_{i=1}^{k}f(a_{i})$
edges, all of them on $Q$. Hence $R$ has length $\ell(R)\geq f(l^{\prime
})+2\sum_{i=1}^{k}f(a_{i})$. Define the broadcast $f^{\prime}$ by%
\[
f^{\prime}(x)=\left\{
\begin{tabular}
[c]{ll}%
$\ell(R)$ & if $x=l^{\prime}$\\
$0$ & if $x\in V(R)-\{l^{\prime}\}$\\
$f(x)$ & otherwise.
\end{tabular}
\right.
\]

Then with respect to $f^{\prime}$, $l^{\prime}$ broadcasts to $b$ but no
farther along $Q$. Together with the fact that each internal vertex of $Q$ has
degree $2$ in $T$, this implies that $f^{\prime}$ is bn-independent. But
\[
\sigma(f^{\prime})=\ell(R)+\sum_{x\in V_{f}^{+}-V(R)}f(x)\geq f(l^{\prime
})+2\sum_{i=1}^{k}f(a_{i})+\sum_{x\in V_{f}^{+}-V(R)}f(x).
\]
Since $k\geq1$,
\[
\sigma(f^{\prime})>f(l^{\prime})+\sum_{i=1}^{k}f(a_{i})+\sum_{x\in V_{f}%
^{+}-V(R)}f(x)=\sigma(f).
\]
This contradicts $f$ being an $\alpha_{\operatorname{bn}}(T)$-broadcast.
Therefore $l^{\prime}$ is the only vertex in $V_{f}^{+}$ on $R$. Since $f$ is
maximal bn-independent, Proposition \ref{prop-max-bn}$(ii)$ implies that $b\in
B_{f}(l^{\prime})$, that is, $f(l^{\prime})=d(l^{\prime},b)$, and the result
follows.~$\blacksquare$

\bigskip

The next lemma states that any $\alpha_{\operatorname{bn}}${-broadcast on a}
tree with at least two branch vertices has at least two broadcasting leaves
belonging to endpaths from different branch vertices, and is useful for
finding an $\alpha_{\operatorname{bn}}${-broadcast on a given tree.}

\begin{lemma}
{\label{Lem_2lb} Any tree $T$ with $b(T)\geq2$ has an }$\alpha
_{\operatorname{bn}}${-broadcast $f$ with two leaves $l,l^{\prime}\in
V_{f}^{+}$ and two distinct branch vertices $w,w^{\prime}\in B(T)$ such that
$l\in L(w)$ and $l^{\prime}\in L(w^{\prime})$.}
\end{lemma}

\noindent\textbf{Proof}.\hspace{0.1in}Let $T$ be a tree with $b(T)\geq2$.
Since $\mathcal{B}(T)$ is a tree of order $b(T)\geq2$, it has at least two
leaves or, equivalently, there are at least two vertices $b_{1},b_{2}\in B(T)$
such that $\operatorname{deg}_{\mathcal{B}(T)}(b_{i})=1$. Thus $|L_{T}%
(b_{i})|\geq2$ for $i=1,2$. Let $g$ be any $\alpha_{\operatorname{bn}}%
(T)$-broadcast. Since $g$ is dominating and leaves only hear leaves (Lemma
\ref{Lem_hear1}), there is at least one broadcasting leaf. We assume that
there is exactly one branch vertex, say $w$, such that $L(w)$ contains all
broadcasting leaves, else $g$ is the required broadcast and our statement is satisfied.

Let $l_{1},...,l_{r}$, $r\geq1$, be the broadcasting leaves in $L(w)$. Any
leaf in $L(w)$ that broadcasts to leaves in $L(T)-L(w)$ overdominates $w$. By
bn-independence, at most one leaf in $L(w)$ overdominates $w$. Hence there is
a unique leaf, say $l_{1}$, in $L(w)$ that dominates $L(T)-L(w)$. Since
$l_{1}$ dominates $L(T)-L(w)$ and since leaves only hear leaves, the only
vertices $l_{1}$ does not dominate lie on $l_{i}-w$ paths, where $2\leq i\leq
r$. By Lemma \ref{Lem_leaf_branch}, these vertices are dominated by the
respective $l_{i}$, $1<i\leq r$. Hence the only broadcasting vertices are
$l_{1},...,l_{r}$. Let $l^{\prime}\in L(T)-L(w)$ be a leaf such that
$d(l^{\prime},w)=\max\{d(l,w):l\in L(T)-L(w)\}$; say $l^{\prime}\in
L(w^{\prime})$, where\ $w^{\prime}\in B(T)-\{w\}$. Since $g(l_{1})\leq
e(l_{1})$, $B_{g}(l_{1})\neq\varnothing$. There are two possible locations for
vertices in $B_{g}(l_{1})$, namely, (a) a leaf in $L(T)-L(w)$, or (b) a vertex
$v\neq w$ on an $l^{\prime\prime}-w$ path, where $l^{\prime\prime}\in L(w)$.
(Possibly, both (a) and (b) hold.) If (a) holds, then, by definition,
$l^{\prime}\in B_{g}(l_{1})$. If (b) holds, then, by Lemma
\ref{Lem_leaf_branch}, $g(l^{\prime\prime})=d(l^{\prime\prime},v)$. (Note that
if $v=l^{\prime\prime}$, then $g(l^{\prime\prime})=0$.) Define a broadcast $f$
by
\[
f(x)=\left\{
\begin{tabular}
[c]{ll}%
$d(x,w)$ & if $x=l^{\prime}$ or $x\in L(w)$\\
$0$ & otherwise.
\end{tabular}
\ \right.
\]

For all leaves $x\in L(w)\cup\{l^{\prime}\}=V_{f}^{+}$, $B_{f}(x)=\{w\}$.
Hence $f$ is bn-independent. If (a) holds, then $f(l_{1})+f(l^{\prime
})=g(l_{1})=d(l_{1},l^{\prime})$, and if (b) holds, then $f(l_{1}%
)+f(l^{\prime\prime})=g(l_{1})+g(l^{\prime\prime})=d(l_{1},l^{\prime\prime})$.
For all $x\in L(w)-\{l_{1},...,l_{r}\}$, $f(x)>g(x)=0$. Suppose $r\geq2$.
Since $l_{1}$ overdominates $w$ in $g$ and $g$ is bn-independent,
$f(l_{i})>g(l_{i})$ for $2\leq i\leq r$. Hence $\sigma(f)\geq\sigma(g)$. Since
$g$ is an $\alpha_{\operatorname{bn}}$-broadcast, so is $f$. Moreover, $f$
satisfies the required conditions (with $l_{1}=l$).~$\blacksquare$

\bigskip

In our final lemma we demonstrate the importance of the branch vertices of a
tree $T$ in determining the broadcast values which may be assigned to its
leaves. By Lemma \ref{Lem_nonleaf1} every tree has an $\alpha
_{\operatorname{bn}}$-broadcast $f$ such that $f(v)=1$ whenever $v\in
V_{f}^{+}-L(T)$. Hence, knowing the broadcast values on the leaves will be
helpful in determining $\alpha_{\operatorname{bn}}(T)$. Combining Lemmas
\ref{Lem_hear1}, \ref{Lem_nonleaf1} and \ref{Lem_branch} allows us to
determine an upper bound for $\alpha_{\operatorname{bn}}(T)$ based on the
order of $T$ and the number and type of branch vertices of $T$.

Suppose a vertex $v$ overdominates a branch vertex $w$ of degree $k$. Once it
has dominated $w$, for the remaining broadcast of $f(v)-d(w,v)$, the broadcast
from $v$ covers up to $k-1$ distinct paths for the same strength of broadcast
required for a single path of this length. Initially, it seems that maximizing
the weight of a broadcast would require dominating each of these $k-1$ paths
with different broadcasting vertices. However, as we will see, sometimes a
maximum weight broadcast is only produced by overdominating branch vertices.

By choosing $\alpha_{\operatorname{bn}}$-broadcasts which minimize the number
of overdominated branch vertices and examining the cost/benefit (for the total
weight) of different types of branch overdomination, Lemma \ref{Lem_branch}
provides restrictions on the way in which a leaf may overdominate a branch
vertex. Informally, Lemma \ref{Lem_branch} states that a leaf $l$ may never
overdominate a branch vertex $b$ by exactly $2$. Either $l$ overdominates a
branch vertex $b$ by exactly $1$ and $b$ has no leaves except possibly $l$, or
$l$ overdominates $b$ by at least $3$ and has exactly one vertex in its
boundary, this vertex being not on a $b^{\prime}-l^{\prime}$ path for any
$b^{\prime}\in B(T)$ and $l^{\prime}\in L(T)$. In addition to $b$, $l$ may
overdominate an unlimited number of branch vertices by $3$ or more as long as
it also overdominates all of their leaves.

We state and prove two claims within the proof; the end of the proof of each
claim is indicated by an open diamond ($\lozenge$).

\begin{lemma}
\label{Lem_branch}Let $T$ be a tree with $b(T)\geq2$. Then there exists an
$\alpha_{\operatorname{bn}}(T)$-broadcast $f$ with the minimum number of
overdominated branch vertices that satisfies the following statement:

For any leaf $l$, let $X$ be the set of all branch vertices overdominated by
$l$. If $X\neq\varnothing$ and $v\in B_{f}(l)$, then $v$ is neither the leaf
nor an internal vertex on any endpath of $T$. Moreover,

\begin{enumerate}
\item[$(i)$] there exists $w\in X$ such that $f(l)=d(l,w)+1$, and either
$L(w)=\{l\}$ and $X=\{w\}$, or $L(w)=\varnothing$ and $f(l)\geq d(l,w^{\prime
})+3$ for all $w^{\prime}\in X-\{w\}$, or

\item[$(ii)$] $f(l)\geq d(l,w)+3$ for all $w\in X$ and $|B_{f}(l)|=1$.
\end{enumerate}
\end{lemma}

\noindent\textbf{Proof}.\hspace{0.1in}Let $f$ be an $\alpha_{\operatorname{bn}%
}(T)$-broadcast for which the number of overdominated branch vertices is a
minimum. Assume that there is a leaf $l$ overdominating a branch vertex. Let
$X=\{w\in B(T):l\ $overdominates$\ w\}\neq\varnothing$ and $A=\{v:v$ is the
leaf or an internal vertex of an endpath of $T\}$. We prove two statements
which we formulate as claims for referencing.\smallskip

\noindent\textbf{Claim \ref{Lem_branch}.1\hspace{0.1in}}\emph{If }$w\in
X$\emph{, then }$f(l)\neq d(l,w)+2$\emph{.}\smallskip

\noindent\textbf{Proof of Claim \ref{Lem_branch}.1}.\hspace{0.1in}Suppose, for
a contradiction, that there exists $w\in X$ such that $f(l)=d(l,w)+2$. Since
$\operatorname{deg}(w)\geq3$, $w$ is adjacent to two vertices $v_{1},v_{2}$
that do not lie on the $l-w$ path. Define the broadcast $g_{1}$ by
$g_{1}(l)=f(l)-2$, $g_{1}(v_{1})=g_{1}(v_{2})=1$ and $g_{1}(u)=f(u)$
otherwise. Then $\sigma(g_{1})=\sigma(f)$. Since $B_{g_{1}}(v_{1})\cap
B_{g_{1}}(v_{2})=\{w\}$, $B_{g_{1}}(l)\cap B_{g_{1}}(v_{i})=\{w\}$ for $i=1,2$
and $N_{g_{1}}(v_{1})\cup N_{g_{1}}(v_{2})\subseteq N_{f}(l)$, $g_{1}$ is
bn-independent. By the maximality of $f$, $g_{1}$ is an $\alpha
_{\operatorname{bn}}(T)$-broadcast. However, since $w$ is no longer
overdominated, $g_{1}$ has fewer overdominated branch vertices than $f$ does,
contrary to the choice of $f$.~$\lozenge$\smallskip

\noindent\textbf{Claim \ref{Lem_branch}.2\hspace{0.1in}}\emph{If }$v\in
B_{f}(l)$\emph{, then }$v$\emph{ is neither the leaf nor an internal vertex of
an endpath of~}$T$\emph{.}\smallskip

\noindent\textbf{Proof of Claim \ref{Lem_branch}.2}.\hspace{0.1in}Suppose, to
the contrary, that $y\in B_{f}(l)\cap A$. Thus there exists $x\in B(T)$ and
$l^{\prime}\in L(T)$ such that $y$ is on an $l^{\prime}-x$ path and
$f(l)=d(l,y)$. Then $l$ overdominates $x$, so $x\in X$. Notice that $l\neq
l^{\prime}$, else $X=\varnothing$. Since $\deg(x)\geq3$, there is a neighbour
$x^{\prime}$ of $x$ such that $x^{\prime}$ does not lie on the $l^{\prime}-x$
path or the $l-x$ path. Since $l$ overdominates $x$, $l$ dominates $x^{\prime
}$. By Lemma \ref{Lem_leaf_branch}, $f(l^{\prime})=d(l^{\prime},y)$. Define
the broadcast $g_{2}$ by $g_{2}(l^{\prime})=d(l^{\prime},x)$, $g_{2}%
(l)=d(l,x)$ and $g_{2}(u)=f(u)$ otherwise. Notice that $N_{g_{2}}(l^{\prime
})\cup N_{g_{2}}(l)\subseteq N_{f}(l)$ and $N_{g_{2}}(l^{\prime})\cap
N_{g_{2}}(l)=\{x\}$. Hence $g_{2}$ is bn-independent. Since $f(l)+f(l^{\prime
})=d(l,l^{\prime})=g_{2}(l)+g_{2}(l^{\prime})$, $g_{2}$ has the same weight as
$f$ and is therefore an $\alpha_{\operatorname{bn}}(T)$-broadcast. However,
$g_{2}$ does not dominate $x^{\prime}$, contradicting Proposition
\ref{prop-max-bn}$(i)$. We conclude that $B_{f}(l)\cap A=\varnothing
$.~$\lozenge$

\smallskip

To show that $(i)$ holds, assume that there exists $w\in X$ such that
$f(l)=d(l,w)+1$. Claim \ref{Lem_branch}.2 implies that $L(w)\subseteq\{l\}$.
If $L(w)=\{l\}$, then $l$ does not overdominate any other branch vertices,
hence $X=\{w\}$. Assume therefore that $L(w)=\varnothing$. Suppose $w^{\prime
}\in X-\{w\}$ and $f(l)=d(l,w^{\prime})+1$. Since $l$ is a leaf,
$d(l,w)=d(l,w^{\prime})\geq2$. Let $v_{1}$ and $v_{2}$ be the neighbours of
$w$ and $w^{\prime}$ on the $l-w$ and $l-w^{\prime}$ paths, respectively,
where possibly $v_{1}=v_{2}$. Create a broadcast $g_{3}$ with $g_{3}%
(l)=f(l)-2$, $g_{3}(w)=g_{3}(w^{\prime})=1$ and $g_{3}(u)=f(u)$ otherwise.
Notice that $N_{g_{3}}(l)\cup N_{g_{3}}(w)\cup N_{g_{3}}(w^{\prime})\subseteq
N_{f}(l)$, $N_{g_{3}}(l)\cap N_{g_{3}}(w)=\{v_{1}\}$, $N_{g_{3}}(l)\cap
N_{g_{3}}(w^{\prime})=\{v_{2}\}$ and either $N_{g_{3}}(w)\cap N_{g_{3}%
}(w^{\prime})=\varnothing$ or $v_{1}=v_{2}$ and $N_{g_{3}}(w)\cap N_{g_{3}%
}(w^{\prime})=\{v_{1}\}$. Hence $g_{3}$ is bn-independent. Notice that
$\sigma(f)=\sigma(g_{3})$ and $g_{3}$ overdominates fewer branch vertices,
contradicting the choice of $f$. Hence $l$ overdominates at most one branch
vertex by exactly one and Claim \ref{Lem_branch}.1 now implies $(i)$.

To show that $(ii)$ holds, assume that $f(l)\geq d(l,w)+3$ for each $w\in X$.
Claim \ref{Lem_branch}.2 implies that $l$ (over)dominates $L(w)$ for each
$w\in X$. Suppose $|B_{f}(l)|\geq2$ and consider two distinct vertices
$v_{1},v_{2}\in B_{f}(l)$. Let $R_{i}$ be the $l-v_{i}$ path for $i=1,2$.
Among all vertices in$\ X\cap V(R_{1})\cap V(R_{2})$, choose $x$ such that
$d(x,l)$ is a maximum. Let $Q_{i}$ be the $x-v_{i}$ subpath of $R_{i}%
,\ i=1,2$. \label{here}By the choice of $x$, $Q_{1}$ and $Q_{2}$ are
internally disjoint. Since $d(l,v_{1})=d(l,v_{2})=f(l)$, $Q_{1}$ and $Q_{2}$
have the same length, say $k$. Since $f(l)\geq d(l,x)+3$, $k\geq3$. Say
$Q_{i}=(x=q_{i,0},q_{i,1},...,q_{i,k}=v_{i}),\ i=1,2$.

\begin{itemize}
\item If $k$ is even, define the broadcast $g_{4}$ by $g_{4}(l)=d(l,x)=f(l)-k$%
, $g_{4}(q_{i,j})=1$ if $j$ is odd and $g_{4}(u)=f(u)$ otherwise. Since
$\bigcup_{i=1}^{2}\bigcup_{j=0}^{k}N_{g_{4}}(q_{i,j})\subseteq N_{f}(l)$,
$g_{4}$ is bn-independent, and since $k$ is even, there are $2(\frac{k}{2})=k$
vertices $q_{i,j}$ in $V_{g_{4}}^{+}$, which implies that $\sigma
(g_{4})=\sigma(f)$. But $g_{4}$ overdominates fewer branch vertices than $f$
does, and we have a contradiction as before.

\item If $k$ is odd, define the broadcast $f^{\prime}$ by $f^{\prime
}(l)=d(l,x)+1$, $f^{\prime}(q_{i,j})=1$ if $j\geq2$ and $j$ is even, and
$f^{\prime}(u)=f(u)$ otherwise. Since $k$ is odd and $k\geq3$, there are
$2(\frac{k-1}{2})=k-1\geq2$ vertices $q_{i,j}$ in $V_{f^{\prime}}^{+}$. As for
$g_{4}$, $f^{\prime}$ is bn-independent and $\sigma(f^{\prime})=\sigma(f)$. If
$f^{\prime}$ is not maximal independent, it can be extended to a
bn-independent broadcast with weight greater than $\sigma(f)$, which is
impossible. Hence $f^{\prime}$ is an $\alpha_{\operatorname{bn}}%
(T)$-broadcast. Either $f^{\prime}$ overdominates fewer branch vertices and
violates the choice of $f$, or $f$ and $f^{\prime}$ overdominate the same
number of branch vertices. In the latter case, since $f^{\prime}(l)=d(l,x)+1$,
we have already shown that $(i)$ holds for $l$ with respect to $f^{\prime}$,
and we consider $f^{\prime}$ instead of~$f$.
\end{itemize}

Hence either $|B_{f}(l)|=1$ and $(ii)$ holds for $f$ and $l$, or there exists
another $\alpha_{\operatorname{bn}}(T)$-broadcast $f^{\prime}$ such that $(i)$
holds for $f^{\prime}$ and $l$.~$\blacksquare$

\bigskip

We restate Lemma \ref{Lem_branch} for broadcasts in which only leaves
broadcast with strength greater than $1$. Corollary \ref{Cor_branch} will be
used in the proof of our upper bound, Theorem \ref{ThmMain}.

\begin{corollary}
\label{Cor_branch}Let $T$ be a tree with $b(T)\geq2$ and $\mathcal{F}^{\prime
}$ the set of all $\alpha_{\operatorname{bn}}(T)$-broadcasts in which only
leaves broadcast with strength greater than $1$. Let $\mathcal{F}$ be the set
of broadcasts in $\mathcal{F}^{\prime}$ with the minimum number of
overdominated broadcasts. Then there exists a broadcast $f\in\mathcal{F}$ that
satisfies the following statement:

For any leaf $l$, let $X$ be the set of all branch vertices overdominated by
$l$. If $X\neq\varnothing$ and $v\in B_{f}(l)$, then $v$ is neither the leaf
nor an internal vertex on any endpath of $T$. Moreover,

\begin{enumerate}
\item[$(i)$] there exists $w\in X$ such that $f(l)=d(l,w)+1$, and either
$L(w)=\{l\}$ and $X=\{w\}$, or $L(w)=\varnothing$ and $f(l)\geq d(l,w^{\prime
})+3$ for all $w^{\prime}\in X-\{w\}$, or

\item[$(ii)$] $f(l)\geq d(l,w)+3$ for all $w\in X$ and $|B_{f}(l)|=1$.
\end{enumerate}
\end{corollary}

\noindent\textbf{Proof}.\hspace{0.1in}By Lemma \ref{Lem_nonleaf1},
$\mathcal{F}^{\prime}\neq\varnothing$ and thus $\mathcal{F}\neq\varnothing$.
Choosing any $f\in\mathcal{F}$, we follow the steps in the proof of Lemma
\ref{Lem_branch} in which we either reduce an existing broadcast strength,
increase the broadcast strength on a leaf, or introduce a new broadcasting
vertex of strength $1$. The resulting broadcast $f^{\prime}$ (which may or may
not be the same as $f$) belongs to $\mathcal{F}$ and also satisfies $(i)$ and
$(ii)$.~$\blacksquare$

\section{Proof of Theorem \ref{ThmMain}}

\label{SecMainThmPf}We are now ready to prove our main result, Theorem
\ref{ThmMain}. The proof proceeds by induction on the number of branch
vertices of the tree and has several cases. We indicate the end of the proof
of each case by a solid diamond ($\blacklozenge$). We consider a tree $T$ and
an $\alpha_{\operatorname{bn}}(T)$-broadcast $f$. Throughout the proof we
consider subtrees $T_{i}$ of $T$ and the restriction of $f$ to $T_{i}$; we
denote the restriction $g_{i}$ of $f$ to $T_{i}$ by $g_{i}=f\upharpoonleft
T_{i}$. In almost all cases these restrictions are broadcasts in the strictest
sense, that is, $g_{i}(x)\leq e_{T_{i}}(x)$ for each $i$; the only exception
occurs in Case 3 of the proof. We restate the theorem for convenience.

\bigskip

\noindent\textbf{Theorem \ref{ThmMain}}\hspace{0.1in}\emph{For any tree }%
$T$\emph{ of order }$n$\emph{,} $\alpha_{\operatorname{bn}}(T)\leq
n-b(T)+\rho(T)$.

\bigskip

\noindent\textbf{Proof}.\hspace{0.1in}By Observation \ref{Ob_Base}, the result
is true for trees with at most two branch vertices. Suppose the result is true
for all trees with fewer than $t$ branch vertices, where $t\geq3$, but false
for at least one tree with $t$ branch vertices. Among all such trees, let $T$
be one of smallest order $n$. By Lemma \ref{Lem_nonleaf1} we may consider
$\alpha_{\operatorname{bn}}(T)$-broadcasts in which only leaves broadcast with
strength exceeding $1$. Let $f$ be such a broadcast in which the number of
overdominated branch vertices is a minimum and such that the statement of
Corollary \ref{Cor_branch} applies. By the choice of $T$, $\alpha
_{\operatorname{bn}}(T)=\sigma(f)>n-b(T)+\rho(T)$.

We first show that every vertex in $V_{f}^{+}$ is a leaf. Then we use
Corollary \ref{Cor_branch} to examine the ways in which endbranch vertices are
dominated. All possibilities lead to contradictions. Since an $\alpha
_{\operatorname{bn}}(T)$-broadcast is dominating, and since
$B_{\operatorname{end}}(T)\neq\varnothing$, the result will follow.

Suppose $f$ has a non-leaf broadcasting vertex. We consider two cases,
depending on its degree.\smallskip

\noindent\textbf{Case A:}\hspace{0.1in}There exists a vertex $b_{0}\in
V_{f}^{+}\cap B(T)$. By the choice of $f$, $f(b_{0})=1$. By Lemma
\ref{Lem_hear1}, $b_{0}$ is not adjacent to a leaf. So either $b_{0}$ is
adjacent to a vertex $v$ of degree $2$, or all neighbours of $b_{0}$ are
branch vertices and $L(b_{0})=\varnothing$.

First assume the former; say $N(v)=\{b_{0},b^{\prime}\}$ and consider the two
subtrees $T_{1},T_{2}$ formed by rejoining $v$ to each component of $T-v$ in
the obvious manner, where $T_{1}$ is the subtree that contains $b_{0}$. For
$i=1,2$, let $g_{i}=f\upharpoonleft T_{i}$ and note that each $g_{i}$ is a
bn-independent broadcast on $T_{i}$. If $b_{0}$ does not belong to $R(T)$,
then $|L_{T}(b_{0})|\geq2$, hence $|L_{T_{1}}(b_{0})|\geq3$ and $b_{0}$ does
not belong to $R(T_{1})$. On the other hand, if $b_{0}$ does belong to $R(T)$,
then $b_{0}$ may or may not belong to $R(T_{1})$. Similarly, if $b^{\prime}$
does not belong to $R(T)$ (possibly $b^{\prime}$ is not even a branch vertex),
then $b^{\prime}$ does not belong to $R(T_{2})$, and if $b^{\prime}$ does
belong to $R(T)$, then $b^{\prime}$ may or may not belong to $R(T_{2})$. Any
other vertex of $T_{i}$ that belongs to $R(T_{i})$ also belongs to $R(T)$.
Therefore $\rho(T)\geq\rho(T_{1})+\rho(T_{2})$. Since $\operatorname{deg}%
_{T}(v)=2$, $b(T)=b(T_{1})+b(T_{2})$. Since $\sigma(f)=\sigma(g_{1}%
)+\sigma(g_{2})$ and, by the assumption on $T$, $\sigma(f)\geq n-b(T)+\rho
(T)+1$, we have%
\begin{align}
\sigma(g_{1})+\sigma(g_{2})  &  =\sigma(f)\geq n-b(T)+\rho(T)+1\nonumber\\
&  \geq|V(T_{1})|+|V(T_{2})|-1-b(T_{1})-b(T_{2})+\rho(T_{1})+\rho
(T_{2})+1\nonumber\\
&  =|V(T_{1})|-b(T_{1})+\rho(T_{1})+|V(T_{2})|-b(T_{2})+\rho(T_{2}).
\label{eqTi}%
\end{align}
But since $b(T_{i})\leq b(T)$ and $|V(T_{i})|<|V(T)|$, the choice of $T$
implies that $\alpha_{\operatorname{bn}}(T_{i})\leq n-b(T_{i})+\rho(T_{i})$,
for $i=1,2$. Hence equality holds throughout (\ref{eqTi}), so
\[
\sigma(g_{i})=\alpha_{\operatorname{bn}}(T_{i})=|V(T_{i})|-b(T_{i})+\rho
(T_{i}),\ i=1,2.
\]
Define the broadcast $g_{1}^{\prime}$ on $T_{1}$ by $g_{1}^{\prime
}(v)=2,\ g_{1}^{\prime}(b_{0})=0$ and $g_{1}^{\prime}(x)=g_{1}(x)$ otherwise.
Since $B_{g_{1}^{\prime}}(v)=B_{g_{1}}(b_{0})-\{v\}$, $g_{1}^{\prime}$ is a
bn-independent broadcast. But $\sigma(g_{1}^{\prime})>\sigma(g_{1}%
)=\alpha_{\operatorname{bn}}(T_{1})$, a contradiction.

We conclude that all neighbours of $b_{0}$ are branch vertices, thus
$L(b_{0})=\varnothing$ and $b_{0}\in R(T)$. Let $b_{1},...,b_{k},\ k\geq3$, be
the neighbours of $b_{0}$ in $T$ and let $T_{1},...,T_{k}$ be the subtrees of
$T$ obtained by rejoining $b_{0}$ to each component of $T-b_{0}$ in the
obvious manner. For $i\in\{1,...,k\}$, let $g_{i}=f\upharpoonleft T_{i}$. As
above, each $g_{i}$ is a bn-independent broadcast on $T_{i}$, hence
$\sigma(g_{i})\leq\alpha_{\operatorname{bn}}(T_{i})$. By the choice of $T$,
$\alpha_{\operatorname{bn}}(T_{i})\leq|V(T_{i})|-b(T_{i})+\rho(T_{i})$, so%
\[
\sum_{i=1}^{k}\sigma(g_{i})\leq\sum_{i=1}^{k}\alpha_{\operatorname{bn}}%
(T_{i})\leq\sum_{i=1}^{k}[|V(T_{i})|-b(T_{i})+\rho(T_{i})].
\]
Hence
\begin{equation}
\sum_{i=1}^{k}\sigma(g_{i})\leq\sum_{i=1}^{k}|V(T_{i})|-\sum_{i=1}^{k}%
b(T_{i})+\sum_{i=1}^{k}\rho(T_{i}). \label{eq_sigma3}%
\end{equation}

Similar to the case of $T_{1}$ and $T_{2}$ above, $R(T_{i})\subseteq
(R(T)-\{b_{0}\})\cap V(T_{i})$, but $b_{0}\notin R(T_{i})$ for each $i$. Hence
$\rho(T)\geq\sum_{i=1}^{k}\rho(T_{i})+1$. By construction, $b(T)=\sum
_{i=1}^{k}b(T_{i})+1$, $n=|V(T)|=\sum_{i=1}^{k}|V(T_{i})|-(k-1)$ and
$\sigma(f)=\sum_{i=1}^{k}\sigma(g_{i})-(k-1)$. By the assumption on $T$ we now
have that
\begin{align*}
\sum_{i=1}^{k}\sigma(g_{i})-k+1  &  =\sigma(f)>n-b(T)+\rho(T)\\
&  \geq\sum_{i=1}^{k}|V(T_{i})|-k+1-\sum_{i=1}^{k}b(T_{i})-1+\sum_{i=1}%
^{k}\rho(T_{i})+1.
\end{align*}
Hence%
\[
\sum_{i=1}^{k}\sigma(g_{i})>\sum_{i=1}^{k}|V(T_{i})|-\sum_{i=1}^{k}%
b(T_{i})+\sum_{i=1}^{k}\rho(T_{i}),
\]
contradicting (\ref{eq_sigma3}). We conclude that no branch vertex of $T$ is a
broadcasting vertex.$~\blacklozenge$\smallskip

\noindent\textbf{Case B:}\hspace{0.1in}There exists a broadcasting vertex $v$
with $\deg(v)=2$. By the choice of $f$, $f(v)=1$. Say $N(v)=\{b_{1},b_{2}\}$
and for $i=1,2$, let $T_{i}$ be the subtree of $T$ obtained by joining $v$ to
$b_{i}$ in $T-v$. Let $g_{i}=f\upharpoonleft T_{i}$. Since $\operatorname{deg}%
(v)=2$, $b(T)=b(T_{1})+b(T_{2})$ and $\rho(T)\geq\rho(T_{1})+\rho(T_{2})$.
Hence, by the induction hypothesis and the choice of $T$,%
\begin{align*}
\sigma(f)  &  =\sigma(g_{1})-\sigma(g_{2})-1\leq\alpha_{\operatorname{bn}%
}(T_{1})+\alpha_{\operatorname{bn}}(T_{2})-1\\
&  \leq|V(T_{1})|+|V(T_{2})|-b(T_{1})-b(T_{2})+\rho(T_{1})+\rho(T_{2})-1\\
&  \leq n-b(T)+\rho(T),
\end{align*}
a contradiction.$~\blacklozenge$\smallskip

Therefore only leaves are broadcasting vertices, that is, $V_{f}^{+}\subseteq
L(T)$. To complete the proof, we show that no branch vertex is overdominated
by exactly $1$ (Case 1). We then consider the way the endbranch vertices are
dominated. By Corollary \ref{Cor_branch}, there are two further ways to
dominate an endbranch vertex $b_{0}$: either $b_{0}$ is dominated but not
overdominated (Case 2), or it is overdominated by $3$ or more (Case 3). We
show that all three cases are impossible.\smallskip

\noindent\textbf{Case 1:\hspace{0.1in}}Suppose that a branch vertex $b_{0}$ is
overdominated by exactly $1$. Since $V_{f}^{+}\subseteq L(T)$, there is a leaf
$l$ such that $f(l)=d(b_{0},l)+1$. By Corollary \ref{Cor_branch}$(i)$,
$L(b_{0})\subseteq\{l\}$. Thus $b_{0}\in R(T)$. Let $b_{1},...,b_{k-1}$,
$k\geq3$, be the neighbours of $b_{0}$ that do not lie on the $b_{0}-l$ path
and note that $\{b_{1},...,b_{k-1}\}\subseteq B_{f}(l)$. For $i\in
\{1,...,k-1\}$, let $T_{i}$ be the subtree of $T$ obtained by joining $b_{0}$
to $b_{i}$ in $T-b_{0}$. Let $T_{k}$ be the tree induced by $N_{f}(l)$. Since
$L(b_{0})\subseteq\{l\}$, no $b_{i}$ is a leaf of $T$; hence each $T_{i}$ is a
proper subtree of $T$. For $i\in\{1,...,k-1\}$, let $g_{i}$ be the broadcast
obtained by first restricting $f$ to $T_{i}$, and then adding $b_{0}$ as
broadcasting vertex of strength $1$. Let $g_{k}=f\upharpoonleft T_{k}$, that
is, $V_{g_{k}}^{+}=\{l\}$ and $g_{k}(l)=f(l)$. Since $B_{g_{i}}(b_{0}%
)\subseteq B_{f}(l)$ for $i\in\{1,...,k-1\}$, each $g_{i}$ is a bn-independent
broadcast on $T_{i}$.

Consider $T_{k}$ and $g_{k}$. Define the broadcast $h$ on $T_{k}$ by
$h(l)=g_{k}(l)-1=f(l)-1$ and $h(b_{i})=1$ for $i\in\{1,...,k-1\}$. Then $h$ is
bn-independent and $\alpha_{\operatorname{bn}}(T_{k})\geq\sigma(h)=\sigma
(g_{k})+k-2$. Since $T_{k}$ is a proper subtree of $T$, a smallest
counterexample,%
\begin{equation}
\sigma(g_{k})\leq\alpha_{\operatorname{bn}}(T_{k})-(k-2)\leq|V(T_{k}%
)|+b(T_{k})-\rho(T_{k})-k+2. \label{eq_gk}%
\end{equation}
By construction, $b(T)=\sum_{i=1}^{k}b(T_{i})$. Since $\operatorname{deg}%
(b_{0})\geq3$ and each $b_{i}$, $i\in\{1,...,k-1\}$, is a leaf in $T_{k}$,
$b_{0}\notin R(T_{k})$; since $b_{0}\in L(T_{i})$ for $i\in\{1,...,k-1\}$,
$b_{0}\notin R(T_{i})$ for $i\in\{1,...,k-1\}$; and since $L_{T}%
(b_{i})\subseteq L_{T_{i}}(b_{i})$ for all $i=1,...,k-1$, each $b_{i}$ that
does not belong to $R(T)$ also does not belong to $R(T_{i})$ for
$i\in\{1,...,k-1\}$. Hence $\rho(T)\geq\sum_{i=1}^{k}\rho(T_{i})+1$. Therefore%
\[
|V(T)|=\sum_{i=1}^{k}|V(T_{i})|-2(k-1),\ \ b(T)=\sum_{i=1}^{k}b(T_{i})\text{
\ and\ \ }\rho(T)\geq\sum_{i=1}^{k}\rho(T_{i})+1.
\]
Since $f(b_{0})=0$ while $g_{i}(b_{0})=1$ for all $i\neq k$, and
$f(x)=g_{i}(x)$ otherwise, $\sigma(f)=\sum_{i=1}^{k}\sigma(g_{i})-(k-1)$.
Therefore%
\[
\sum_{i=1}^{k}\sigma(g_{i})-(k-1)=\sigma(f)>n-b(T)+\rho(T)\geq\sum_{i=1}%
^{k}|V(T_{i})|-2(k-1)-\sum_{i=1}^{k}b(T_{i})+\sum_{i=1}^{k}\rho(T_{i})+1
\]
and so%
\begin{equation}
\sum_{i=1}^{k}\sigma(g_{i})>\sum_{i=1}^{k}|V(T_{i})|-\sum_{i=1}^{k}%
b(T_{i})+\sum_{i=1}^{k}\rho(T_{i})-k+2. \label{eq_sigma6}%
\end{equation}
Since each $T_{i}$ is a proper subtree of $T$, $\sigma(g_{i})\leq
\alpha_{\operatorname{bn}}(T_{i})\leq|V(T_{i})|+b(T_{i})-\rho(T_{i})$ for
$i\in\{1,...,k-1\}$, and by (\ref{eq_gk}), $\sigma(g_{k})\leq\alpha
_{\operatorname{bn}}(T_{k})-k+2$. Hence
\[
\sum_{i=1}^{k}\sigma(g_{i})\leq\sum_{i=1}^{k}|V(T_{i})|-\sum_{i=1}^{k}%
b(T_{i})+\sum_{i=1}^{k}\rho(T_{i})-k+2,
\]
which contradicts (\ref{eq_sigma6}). Hence no branch vertex is overdominated
by exactly $1$.$~\blacklozenge$\smallskip

We now focus on the possible ways to dominate an endbranch vertex.\smallskip

\noindent\textbf{Case 2:\hspace{0.1in}}Suppose there exists a vertex $b_{0}\in
B_{\operatorname{end}}(T)$ that is dominated, but not overdominated, by a leaf
$l$. Then $f(l)=d(b_{0},l)$ and $b_{0}\in B_{f}(l)$. By Lemma
\ref{Lem_leaf_branch}, $f(l^{\prime})=d(b_{0},l^{\prime})$ for each
$l^{\prime}\in L(b_{0})$. It is possible that the only leaves dominating
$b_{0}$ are in $L(b_{0})$. Let $v_{1}$ be the neighbour of $b_{0}$ that does
not lie on a $b_{0}-l^{\prime}$ path for any $l^{\prime}\in L(b_{0})$. Since
$b_{0}$ is not overdominated, no vertex in $L(b_{0})$ dominates $v_{1}$. Let
$T_{0}$ be the subtree of $T-b_{0}v_{1}$ that contains $b_{0}$ and let $T_{1}$
be the subtree obtained by joining $b_{0}$ to $v_{1}$ in the subtree of
$T-b_{0}v_{1}$ that contains $v_{1}$. Then $T_{0}$ is a path or a generalized
spider. By Proposition \ref{Prop=}, $\alpha_{\operatorname{bn}}(T_{0}%
)=|V(T_{0})|-1\geq|V(T_{0})|-b(T_{0})+\rho(T_{0})$. Also, $b(T_{1})=b(T)-1$
and, since $b_{0}$ is an endbranch vertex, $\rho(T_{1})\leq\rho(T)$. We
consider two subcases, depending on whether some vertex in $T_{1}$ dominates
$b_{0}$ or not.\smallskip

\noindent\textbf{Subcase 2.1:\hspace{0.1in}}No vertex in $T_{1}$ dominates
$b_{0}$. Let $f_{1}=f\upharpoonleft T_{1}$ and define $g_{1}$ by $g_{1}%
(b_{0})=1$ and $g_{1}(x)=f_{1}(x)$ otherwise. Let $g_{0}=f\upharpoonleft
T_{0}$. Then%
\begin{align*}
\sigma(g_{0})+\sigma(g_{1})-1  &  =\sigma(f)>n-b(T)+\rho(T)\\
&  \geq|V(T_{0})|+|V(T_{1})|-1-(b(T_{1})+1)+\rho(T_{1}),
\end{align*}
hence
\begin{equation}
\sigma(g_{0})+\sigma(g_{1})>|V(T_{0})|+|V(T_{1})|-b(T_{1})+\rho(T_{1})-1.
\label{eq_sigma5}%
\end{equation}
But $T_{1}$ is a proper subtree of $T$, hence, by the choice of $T$,
$\sigma(g_{1})\leq\alpha_{\operatorname{bn}}(T_{1})\leq|V(T_{1})|-b(T_{1}%
)+\rho(T_{1})$. Combined with Proposition \ref{Prop=} this gives%
\[
\sigma(g_{0})+\sigma(g_{1})\leq|V(T_{0})|+|V(T_{1})|-b(T_{1})+\rho(T_{1})-1,
\]
which contradicts (\ref{eq_sigma5}).~$\lozenge$\smallskip

\noindent\textbf{Subcase 2.2:\hspace{0.1in}}Some vertex $y$ in $T_{1}$
dominates $b_{0}$. Since $V_{f}^{+}\subseteq L(T)$ (as proved in Cases A and
B), $y$ is a leaf. By the assumption for Case 2, $v_{1}$ does not lie on a
$b_{0}-l^{\prime}$ path for $l^{\prime}\in L(b_{0})$, and $y$ does not
overdominate $b_{0}$. Hence $b_{0}\in B_{f}(y)$ and $f(y)=d(b_{0}%
,y)=d(v_{1},y)+1$; moreover, the $v_{1}-y$ path contains a branch vertex. Let
$c$ be the branch vertex on this path nearest to $v_{1}$. As shown in Case 1,
no leaf overdominates a branch vertex by exactly $1$, hence $v_{1}$ is not a
branch vertex and therefore $c\neq v_{1}$. By Corollary \ref{Cor_branch},
$f(y)\geq d(c,y)+3$ and $B_{f}(y)$ consists of a single non-leaf vertex. But
$b_{0}\in B_{f}(y)$, hence $B_{f}(y)=\{b_{0}\}$. This, however, implies that
$y$ dominates all of $T_{1}$, otherwise $B_{f}(y)$ would contain another
vertex. For $i=0,1$, let $g_{i}=f\upharpoonleft T_{i}$. We now have that
$\sigma(g_{1})=f(y)=e_{T_{1}}(y)\leq\operatorname{diam}(T_{1})\leq
|V(T_{1})|-b(T_{1})-1$ (by Observation \ref{Obs_diam}), hence%
\[
\sigma(g_{0})+\sigma(g_{1})\leq|V(T_{0})|-1+|V(T_{1})|-b(T_{1})-1=|V(T_{0}%
)|+|V(T_{1})|-b(T_{1})-2.
\]
However, since $\rho(T_{0})=0$, $\rho(T_{1})\leq\rho(T)$ and $b(T)=b(T_{1}%
)+1$,%
\begin{align*}
\sigma(g_{0})+\sigma(g_{1})  &  =\sigma(f)>n-b(T)+\rho(T)\\
&  \geq|V(T_{0})|+|V(T_{1})|-1-(b(T_{1})+1)+\rho(T_{1})\\
&  =|V(T_{0})|+|V(T_{1})|-b(T_{1})+\rho(T_{1})-2.
\end{align*}
This contradiction concludes the proof of Subcase 2.2 and thus the proof of
Case 2. Hence no $v\in B_{\operatorname{end}}(T)$ is dominated without being
overdominated.$~\blacklozenge$\smallskip

We have shown that each end-branch vertex that is dominated by a leaf $l$ is
overdominated by more than $1$ by $l$. By Corollary \ref{Cor_branch}, only one
case remains to be considered.\smallskip

\noindent\textbf{Case 3:\hspace{0.1in}}A vertex $b_{0}\in
B_{\operatorname{end}}(T)$ is overdominated by a leaf $l$ and $f(l)\geq
d(l,b)+3$. Let $L(b_{0})=\{l_{1},...,l_{k}\}$. Since $b_{0}$ is an endbranch
vertex, $k\geq2$. There are two subcases: $l\in L(b_{0})$ or $l\notin
L(b_{0})$. In either case, by Corollary \ref{Cor_branch}, $B_{f}(l)=\{v\}$ for
some vertex $v$ which is neither the leaf nor an internal vertex on any
endpath of $T$. Let $v^{\prime}$ be the vertex on the $v-l$ path such that
$d(v,v^{\prime})=2$. Then $l$ overdominates $v^{\prime}$ by exactly $2$. By
Corollary \ref{Cor_branch}, $v^{\prime}$ is not a branch vertex, hence
$\operatorname{deg}(v^{\prime})=2$. Form two subtrees of $T$ by reconnecting
$v^{\prime}$ to each component of $T-\{v^{\prime}\}$ in the obvious way. Let
$T_{1}$ be the tree which contains $l$ and $T_{2}$ the other tree. Let
$f_{i}=f\upharpoonleft T_{i}$ for $i=1,2$. Notice that the $v^{\prime}-v$ path
is not $f_{2}$-dominated. Extend $f_{2}$ by creating a broadcast
$f_{2}^{\prime}$ on $T_{2}$ with $f_{2}^{\prime}(v^{\prime})=2$ and
$f_{2}^{\prime}(x)=f_{2}(x)$ otherwise. Since $f$ is bn-independent and
$N_{f_{2}^{\prime}}(v^{\prime})\subset N_{f}(l)$, $f_{2}^{\prime}$ is
bn-independent and
\begin{equation}
\sigma(f_{2})+2=\sigma(f_{2}^{\prime})\leq\alpha_{\operatorname{bn}}(T_{2}).
\label{eq_min2}%
\end{equation}

\noindent\textbf{Case 3.1:\hspace{0.1in}}Suppose $l\in L(b_{0})$. Without loss
of generality, say $l=l_{1}$. Since $l_{1}$ overdominates $b_{0}$ and
$B_{f}(l_{1})=\{v\}$, $l_{1}$ dominates $T_{1}$. (Since $f_{1}(l)=e_{T_{1}%
}(l)+2$, $f_{1}$ is not a broadcast on $T_{1}$). Define a new broadcast
$f_{1}^{\prime}$ on $T_{1}$ with $f_{1}^{\prime}(l_{i})=d(l_{i},b_{0})$ for
$i=1,...,k$, $f_{1}^{\prime}(v^{\prime})=d(v^{\prime},b_{0})$ and
$f_{1}^{\prime}(x)=0$ otherwise. For all $1\leq i,j\leq k$, $i\neq j$,
$N_{f_{1}^{\prime}}(l_{i})\cap N_{f_{1}^{\prime}}(l_{j})=B_{f_{1}^{\prime}%
}(l_{i})\cap B_{f_{1}^{\prime}}(l_{j})=\{b_{0}\}$ and $N_{f_{1}^{\prime}%
}(l_{i})\cap N_{f_{1}^{\prime}}(v^{\prime})=B_{f_{1}^{\prime}}(l_{i})\cap
B_{f_{1}^{\prime}}(v^{\prime})=\{b_{0}\}$. Hence $f_{1}^{\prime}$ is a
bn-independent broadcast. Since $d(v,v^{\prime})=2$, it follows that
$\sigma(f_{1})=f_{1}(l_{1})=d(l_{1},v^{\prime})+2=f_{1}^{\prime}(l_{1}%
)+f_{1}(v^{\prime})+2$. Since we also have that $f_{1}^{\prime}(l_{i}%
)=d(l_{i},b_{0})\geq1$ for all $i=2,...,k$, we deduce that $\sigma
(f_{1}^{\prime})\geq\sigma(f_{1})-2+k-1$. Since $k\geq2$,
\[
\sigma(f_{1})\leq\sigma(f_{1}^{\prime})+3-k\leq\sigma(f_{1}^{\prime}%
)+1\leq\alpha_{\operatorname{bn}}(T_{1})+1.
\]
Further, both graphs have fewer vertices than $T$ and $b(T_{1}),b(T_{2}%
)<b(T)$. Hence, by the induction hypothesis, $\sigma(f_{1})\leq V(T_{1}%
)-b(T_{1})+\rho(T_{1})+1$ and $\sigma(f_{2})+2=\sigma(f_{2}^{\prime})\leq
V(T_{2})-b(T_{2})+\rho(T_{2})$. It follows that
\[
\sigma(f_{1})+\sigma(f_{2})\leq|V(T_{1})|+|V(T_{2})|-b(T_{1})-b(T_{2}%
)+\rho(T_{1})+\rho(T_{2})-1.
\]
Since $L(T)\subseteq L(T_{1})\cup L(T_{2})$, $\rho(T)\geq\rho(T_{1}%
)+\rho(T_{2})$, and since $\operatorname{deg}(v)=2$, $b(T_{1})+b(T_{2})=b(T)$.
By construction, $|V(T_{1})|+|V(T_{2})|=n+1$. Hence%

\[
\sigma(f)=\sigma(f_{1})+\sigma(f_{2})\leq n+1-b(T)+\rho(T)-1,
\]
contradicting our choice of $f$.\smallskip

\noindent\textbf{Case 3.2:\hspace{0.1in}}Suppose $l\in L(b)$ where $b\neq
b_{0}$. Since $l$ overdominates $b_{0}$ and $B_{f}(l)=\{v\}$, where $v$ is
neither the leaf nor an internal vertex on an endpath, $l$ overdominates
$L(b_{0})$. Define the broadcast $h_{1}$ on $T_{1}$ by $h_{1}%
(l)=d(l,b),\ h_{1}(l_{1})=d(l_{1},b),\ h_{1}(v^{\prime})=d(v^{\prime},b)$ and
$h(x)=0$ otherwise. Note that, for any two distinct vertices $x,y\in V_{h_{1}%
}^{+},\ N_{h_{1}}(x)\cap N_{h_{1}}(y)=\{b\}$. Hence $h_{1}$ is bn-independent
and $\sigma(h_{1})\leq\alpha_{\operatorname{bn}}(T_{1})$. Since $b_{0}\neq b$,
$h_{1}(l_{1})\geq2$. Also, $h_{1}(l)+h_{1}(v^{\prime})+2=f_{1}(l)$, so that
$\sigma(h_{1})\geq\sigma(f_{1})$. By the induction hypothesis, and because
$\sigma(f_{2})\leq\alpha_{\operatorname{bn}}(T_{2})-2$ by (\ref{eq_min2}),
\[
\sigma(f_{1})+\sigma(f_{2})\leq\sigma(h_{1})+\sigma(f_{2})\leq|V(T_{1}%
)|+|V(T_{2})|-b(T_{1})-b(T_{2})+\rho(T_{1})+\rho(T_{2})-2.
\]
As before, $\rho(T)\geq\rho(T_{1})+\rho(T_{2})$ and $b(T_{1})+b(T_{2})=b(T)$.
Hence%
\[
\sigma(f)=\sigma(f_{1})+\sigma(f_{2})\leq n+1-b(T)+\rho(T)-2,
\]
contradicting our choice of $f$.~$\blacklozenge$\smallskip

This final contradiction shows that the endbranch vertices of $T$ are not
dominated. Since all maximal bn-independent broadcasts are dominating, the
theorem follows.~$\blacksquare$

\bigskip%

\begin{figure}[ptb]%
\centering
\includegraphics[
height=1.7815in,
width=3.0139in
]%
{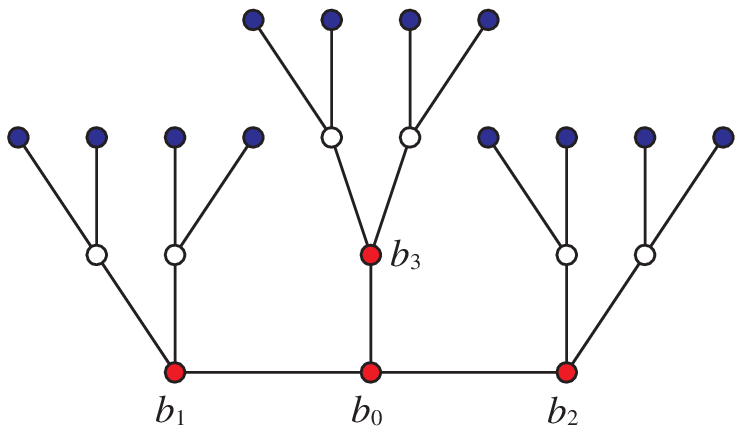}%
\caption{A tree $T$ with $\alpha_{\operatorname{bn}}(T)=15<b-b(T)+\rho(T)=16$}%
\label{Fig_less than ub}%
\end{figure}
Corollary \ref{Cor_Cat} below illustrates that the bound in Theorem
\ref{ThmMain} is sharp. There also exist trees for which strict inequality
holds. Consider the tree $T$ in Figure \ref{Fig_less than ub}, for example.
The vertices in $R(T)$ are shown in red; $\rho(T)=4,\ b(T)=10$ and
$n=|V(T)|=22$, hence the bound is $\alpha_{\operatorname{bn}}(T)\leq16$.
However, the blue leaves together with $\{b_{1},b_{2},b_{3}\}$ form an
$\alpha(T)$-set of cardinality $15$, hence its characteristic function is a
bn-independent broadcast. By considering all possible bn-independent
broadcasts where leaves only hear leaves, it can be shown that indeed
$\alpha_{\operatorname{bn}}(T)=15$.

\section{A brief look at caterpillars}

\label{Sec_Cat}A \emph{caterpillar} of length $k\geq0$ is a tree such that
removing all leaves produces a path of length $k$, called the \emph{spine}. A
vertex on the spine is called a \emph{spine} \emph{vertex}. A caterpillar with
exactly one spine vertex is a star, and one with exactly one branch vertex is
a generalized spider, hence we consider caterpillars with two or more branch
vertices. We think of a caterpillar as drawn with the spine on a horizontal
line, so that we can refer to its leftmost or rightmost branch vertex/spine
vertex. The notation used in what follows is defined in Section
\ref{SecTreeDefs}. If $v$ is a branch vertex of a caterpillar $T$, then $v$ is
a stem, hence $B_{0}(T)=\varnothing$ and $R(T)=B_{1}(T)$.

\begin{figure}[ptb]%
\centering
\includegraphics[
height=1.081in,
width=3.3512in
]%
{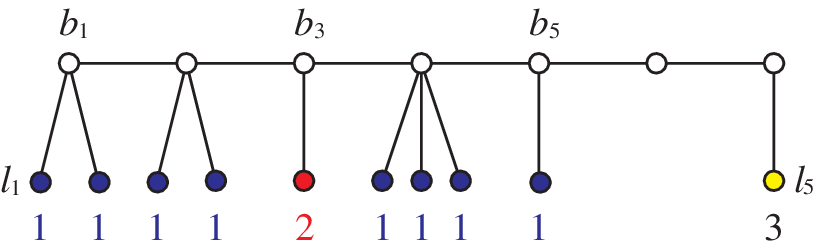}%
\caption{An $\alpha_{\operatorname{bn}}(T)$-broadcast on a caterpillar $T$
with $W_{\operatorname{int}}(T)=\varnothing$ and $B_{1}(T)=R(T)=\{b_{3}\}$}%
\label{Fig_Cat}%
\end{figure}

We consider caterpillars $T$ such that $W_{\operatorname{int}}(T)=\varnothing
$, i.e., there are no vertices of degree $2$ between the leftmost and
rightmost branch vertices, and $B_{1}(T)$ is either empty or an independent
set. (See Figure \ref{Fig_Cat} for an example.) For such a caterpillar $T$,
let $b_{1},...,b_{k}$ be the branch vertices of $T$, labelled from left to
right on the spine. Note that $b_{1},b_{k}\in B_{\geq2}(T)$, i.e.,
$b_{1},b_{k}\notin R(T)$. Let $l_{1}\in L(b_{1})$ and $l_{k}\in L(b_{k})$ be
leaves farthest from $b_{1}$ and $b_{k}$, respectively.%

Define the broadcast $f$ on $T$ by $f(l_{1})=d(l_{1},b_{1}),\ f(l_{k}%
)=d(l_{k},b_{k}),\ f(l)=2$ if $l\in L(b)$ where $b\in B_{1}(T)$, $f(l)=1$ if
$l$ is any other leaf, and $f(v)=0$ for all other vertices $v$. Suppose $l\in
L(b)$. If $b\in B_{1}(T)$, then $B_{f}(l)$ consists of the two branch vertices
on either side of $b$. In all other cases, $B_{f}(l)=\{b\}$. Since $B_{1}(T)$
is independent, $f$ overlaps only in boundaries, hence $f$ is bn-independent.
Since $V_{f}^{+}=L(T)$ and $f(l_{1})$ equals the number of vertices on the
$l_{1}-b_{1}$ path not counting $b_{1}$, and similarly for $l_{k}$, it follows
that $\alpha_{\operatorname{bn}}(T)\geq\sigma(f)=|V(T)|-|B(T)|+|B_{1}%
(T)|=|V(T)|-b(T)+\rho(T)$. Hence we have the following corollary to Theorem
\ref{ThmMain}, which gives a class of trees for which equality holds in the bound.

\begin{corollary}
\label{Cor_Cat}If $T$ is a caterpillar whose branch vertices induce a path
$P=(b_{1},...,b_{k})$, and $R(T)$ (i.e., the branch vertices among
$b_{2},...,b_{k-1}$ that are adjacent to exactly one leaf) is either empty or
an independent set, then $\alpha_{\operatorname{bn}}(T)=|V(T)|-b(T)+\rho(T)$.
\end{corollary}

\section{Open problems}

\label{Sec_Open}For the tree $T$ in Figure \ref{Fig_less than ub},
$\alpha_{\operatorname{bn}}(T)=15<|V(T)|-b(T)+\rho(T)=16$. The subgraph of $T$
induced by $R(T)$ is $K_{1,3}$, and $\alpha(K_{1,3})=3$. Note that
$|V(T)|-b(T)+\alpha(T[R(T)]=15=\alpha_{\operatorname{bn}}(T)$. This raises the
following question.

\begin{question}
Can the upper bound in Theorem \ref{ThmMain} be improved to $\alpha
_{\operatorname{bn}}(T)\leq|V(T)|-b(T)+\alpha(T[R(T)]$?
\end{question}

\begin{problem}
Characterize trees $T$ such that $\alpha_{\operatorname{bn}}%
(T)=|V(T)|-b(T)+\rho(T)$.
\end{problem}

Ahmane et al.\ \cite{ABS} determined $\alpha_{h}(T)$ for caterpillars $T$ for
which $W(T)$ is an independent set. (Here we regard the empty set as being
independent.) It appears that their proof also works if we only require that
$W_{\operatorname{int}}(T)$ is independent. What is interesting in view of the
fact that $\alpha_{h}/\alpha_{\operatorname{bn}}<2$, the ratio being
asymptotically best possible even for trees \cite{MN}, is that $\alpha
_{\operatorname{bn}}(T)=\alpha_{h}(T)$ for some classes of caterpillars. They
showed (see \cite[Corollary 18]{ABS}) that if $W(T)=\varnothing$ (or perhaps
we only need that $W_{\operatorname{int}}(T)=\varnothing$) and at least one
stem is adjacent to three or more leaves, then $\alpha_{h}(T)=|V(T)|-b(T)+\rho
(T)$, and (in Corollary 20) that if $W(T)$ is independent (or perhaps if
$W_{\operatorname{int}}(T)$ is independent) and all stems (except perhaps the
leftmost and rightmost stems) are adjacent to three or more leaves (thus
$\rho(T)=0$), then $\alpha_{h}(T)=|V(T)|-b(T)$. In either case, $\alpha
_{h}(T)=\alpha_{\operatorname{bn}}(T)$.

\begin{problem}
Characterize caterpillars $T$ such that $\alpha_{\operatorname{bn}}%
(T)=\alpha_{h}(T)$.
\end{problem}

\begin{problem}
Determine $\alpha_{\operatorname{bn}}(T)$ for all caterpillars $T$.
\end{problem}

\begin{problem}
Characterize trees $T$ such that $\alpha_{\operatorname{bn}}(T)=\alpha_{h}(T)$.
\end{problem}

\label{refs}

\end{document}